\documentclass[10pt]{amsart}
\usepackage{amsmath, amssymb, latexsym, enumerate, graphicx}

\newtheorem{thm}[equation]{Theorem}
\newtheorem{pro}[equation]{Proposition}

\newtheorem{cor}[equation]{Corollary}
\newtheorem{lem}[equation]{Lemma}

\makeatletter\@addtoreset{equation}{section}\makeatother

\theoremstyle{definition}
\newtheorem{exa}[equation]{Example}

\newcommand{\spandsp}{\qquad\text{and}\qquad}

\newcommand{\co}{\colon}
\newcommand{\rta}{\rightarrow}
\renewcommand{\o}{\otimes}
\renewcommand{\=}{\,=\,}
\newcommand{\+}{\,+\,}
\newcommand{\isom}{\cong}
\newcommand{\pr}[2]{\ensuremath{\langle {#1,#2}\rangle}}
\newcommand{\sub}[1]{_{\text{${\scriptscriptstyle #1}$}}}

\newcommand{\rkm}[1]{\ensuremath{\rho_{\scriptscriptstyle M}(#1)}}
\newcommand{\rkn}[1]{\ensuremath{\rho_{\scriptscriptstyle N}(#1)}}
\newcommand{\rkl}[1]{\ensuremath{\rho_{\scriptscriptstyle L}(#1)}}
\newcommand{\nln}[1]{\ensuremath{\nu_{\scriptscriptstyle N}(#1)}}
\newcommand{\cl}[1]{\ensuremath{e_{\scriptscriptstyle #1}}}

\newcommand{\frp}{\mathbin{\Box}}
\newcommand{\cc}{\ensuremath{\mathcal C}}
\newcommand{\ci}{\ensuremath{\mathcal I}}
\newcommand{\cm}{\ensuremath{\mathcal M}}
\newcommand{\cp}{\ensuremath{\mathcal P}}
\newcommand{\ccr}{\ensuremath{\mathcal R}}
\newcommand{\cu}{\ensuremath{\mathcal U}}
\newcommand{\cw}{\ensuremath{\mathcal W}}
\newcommand{\cws}{\ensuremath{\mathcal W_{\scriptscriptstyle S}}}
\newcommand{\cwm}{\ensuremath{\mathcal W_{\scriptscriptstyle M}}}
\newcommand{\crm}{\ensuremath{\mathcal R_{\scriptscriptstyle M}}}
\newcommand{\crs}{\ensuremath{\mathcal R_{\scriptscriptstyle S}}}

\newcommand{\cim}{\ensuremath{\mathcal I_{\scriptscriptstyle M}}}
\newcommand{\cis}{\ensuremath{\mathcal I_{\scriptscriptstyle S}}}
\newcommand{\cmi}{\ensuremath{\mathcal M_{\scriptscriptstyle I}}}

\newcommand{\ic}[1]{\ensuremath{\text{\rm IC}(#1)}}
\newcommand{\frmod}[1]{\ensuremath{K\{{#1}\}}}      % Free K-module
\newcommand{\dual}{\ensuremath{^*}}
\newcommand{\un}[2]{\ensuremath{U_{#1,#2}}} % The uniform matroid

\begin{document}
\title[Primitive elements in the matroid-minor Hopf algebra]
{Primitive elements in the matroid-minor Hopf algebra}
\author{Henry Crapo}
\author{William Schmitt}
\thanks{Schmitt partially supported by NSA grant 02G-134}
\email{crapo@ehess.fr and wschmitt@gwu.edu}
\bibliographystyle{amsplain}

\begin{abstract}
  We introduce the {\it matroid-minor} coalgebra $C$, which has
  labeled matroids as distinguished basis, and coproduct given by
  splitting a matroid into a submatroid and complementary contraction
  all possible ways.  We introduce two new bases for $C$; the first of
  these is is related to the distinguished basis by M\"obius inversion
  over the rank-preserving weak order on matroids, the second by
  M\"obius inversion over the suborder excluding matroids that are
  irreducible with respect to the free product operation.  We show
  that the subset of each of these bases corresponding to the set of
  irreducible matroids is a basis for the subspace of primitive
  elements of $C$.  Projecting $C$ onto the matroid-minor Hopf algebra
  $H$, we obtain bases for the subspace of primitive elements of $H$.
\end{abstract}

\maketitle 

\section{Introduction}

Many of the Hopf algebras now of central importance in algebraic
combinatorics share certain striking features, suggesting the
existence of a natural, yet-to-be-identified, class of combinatorial
Hopf algebras. These Hopf algebras are graded and cofree, each has a
canonical basis consisting of, or indexed by, a family of (equivalence
classes of) combinatorial objects that is equipped with a natural
partial ordering, and in each case the algebraic structure is most
clearly understood through the introduction of a second basis, related
to the canonical one by M\"obius inversion over the partial ordering.
Indeed, in a number of key examples cofreeness becomes apparent once
the coproduct is expressed in terms of the second basis, and this
basis contains, as an easily recognizable subset, a basis for the
subspace of primitive elements.  The significance of primitive
elements in this context was established by Loday and Ronco in
\cite{loro:osc}, where they proved a Milnor-Moore type theorem
characterizing cofree Hopf algebras in terms of their primitive
elements.

Examples of such Hopf algebras include the algebra of quasisymmetric
functions, introduced by Gessel in \cite{ge:mpp}, the Hopf algebra
structure of which was determined by Malvenuto in \cite{ma:pcf}, the
Malvenuto-Reutenauer Hopf algebra of permutations \cite{ma:pcf,
  mare:dqf}, the Loday-Ronco Hopf algebra of planar binary trees
\cite{loro:hap}, and, most recently, the Hopf algebra of uniform block
partitions, due to Aguiar and Orellana \cite{agor:hau}.  The canonical
basis for quasisymmetric functions is indexed by compositions of
nonnegative integers, whose natural partial ordering, given by
refinement, is a disjoint union of Boolean algebras.  In
\cite{agso:smr}, Aguiar and Sottile use the weak order on the
symmetric groups to elucidate the structure of the
Malvenuto-Reutenauer Hopf algebra, and in \cite{agso:slr} they use the
Tamari order on planar binary trees in an analogous fashion to study
the Loday-Ronco Hopf algebra.  Furthermore, through the use of Galois
connections between each pair of the aforementioned partial orders,
they exhibit the myriad relationships among these Hopf algebras in a
completely unified manner.

In this article, we use similar techniques to study a seemingly
unrelated Hopf algebra, based on matroids.
 The {\it matroid-minor Hopf algebra}, introduced in
\cite{sc:iha}, has as canonical basis the set of all isomorphism
classes of matroids, with product induced by the direct sum operation,
and coproduct of the isomorphism class \cl M\ of a matroid $M=M(S)$
given by $\sum\sub{A\subseteq S}\cl{M|A}\o\cl{M/A}$, where $M|A$ is
the submatroid obtained by restriction to $A$ and $M/A$ is the
complementary contraction. The current authors showed in
\cite{crsc:uft} that this Hopf algebra is cofree, and that 
its subspace of primitive elements has a basis
indexed by those isomorphism classes of matroids that are irreducible
with respect to the {\it free product} operation \cite{crsc:fpm}.
We approach the matroid-minor Hopf algebra here by first lifting to
the matroid-minor {\it coalgebra}, which has as canonical basis the
set of all {\it labeled} matroids whose underlying sets are subsets of
some given infinite set.  The coproduct of a matroid $M(S)$ in the
matroid-minor coalgebra is given by $\sum\sub{A\subseteq S}M|A\o M/A$,
so the natural projection, taking a matroid to its isomorphism class,
is a coalgebra map.  The set of labeled matroids is partially ordered
by the (rank-preserving) {\it weak order}, under which $M\geq N$ means
that $M$ and $N$ have the same underlying set, and each basis for $N$
is also a basis for $M$.  We introduce two new bases for the
matroid-minor coalgebra, both related to the canonical basis by
M\"obius inversion; the first over the full weak order and the second
over a suborder that excludes matroids which are irreducible with
respect to free product.  The subset of each of these new bases
corresponding to the irreducible matroids is a basis for the subspace of
primitive elements.  Applying the projection map, we obtain bases for
the subspace of primitive elements of the matroid-minor Hopf algebra,
one of these is the basis previously identified in \cite{crsc:uft}.

\section{Posets and M\"obius functions}

In this section we gather together for later use some basic facts
about partially ordered sets and their M\"obius functions.  Of the
four results given here, two (Theorems \ref{thm:closure} and
\ref{thm:hall}) are classical, and stated without proof, one
(Proposition \ref{pro:mobius}) is apparently new, and one (Proposition
\ref{pro:Pzero}) is trivial, but nonetheless useful to have on hand.

We assume  that all partially ordered sets, or {\it posets}, for
short, are locally finite, that is, given $x\leq z$ in a poset $P$,
the interval $[x,z] = \{y\in P\co x\leq y\leq z\}$ is finite.
The {\it M\"obius function} $\mu=\mu\sub P$ of a poset $P$ is the
integer-valued function having the set of intervals in $P$ as domain, 
defined by $\mu
(x,x) = 1$, for all $x\in P$, and
$$
\mu (x,z)\=\, -\sum_{x\leq y< z}\mu (x,y)\=\, 
-\sum_{x< y\leq z}\mu (y,z),
$$
for all $x<z$ in $P$.

A {\it closure operator} on a poset $P$ is an idempotent, order-preserving
map $\varphi\co P\rta P$ such that
$x\leq\varphi (x)$, for all $x\in P$. Given a closure operator $\varphi$
on $P$, we write $P_\varphi$ for the subposet $\text{\rm im\,} \varphi =
\{x\in P\co x=\varphi(x)\}$ of {\it closed elements\/} of $P$.
An essential ingredient in the proofs of our main results,
Theorems \ref{thm:primbasis1} and \ref{thm:primbasis2}, is
the following well-known theorem, due to Rota \cite{ro:fct1}, that
expresses the M\"obius function of $P_\varphi$ in terms of that of $P$.
\begin{thm}\label{thm:closure}
   If $\varphi$ is a closure operator on a poset $P$, then
for all $a\leq b$ in $P$,
$$
    \sum_{x\colon\varphi (x) =b}\mu (a,x)\=
\begin{cases}
      \mu_\varphi (a,b) & \text{if $a,b\in P_\varphi$}, \\
  0 & \text{otherwise},
\end{cases}
$$
where $\mu$ and $\mu_\varphi$ denote the M\"obius functions of $P$ and 
$P_\varphi$, respectively. 
\end{thm}

Given $x\leq z$ in a poset $P$, a {\it chain} from $x$ to $z$ is
a sequence $C=(x_0,\ldots,x_k)$ of elements of $P$ such that
$x=x_0<\cdots<x_k=z$.  The {\it length} $\ell (C)$ of a chain $C$
is $|C|-1$.  The following theorem, due to Philip Hall \cite{ha:ctg},
provides an alternative definition of M\"obius function that
allows us to give a short proof of Proposition \ref{pro:mobius}
below.
\begin{thm}\label{thm:hall}
  If $x\leq z$ in a poset $P$, then $\mu\sub P(x,z)=\sum (-1)^{\ell (C)}$,
where the sum is over all chains $C$ from $x$ to $z$ in $P$.
\end{thm}

Suppose that $P$ is a poset and $Q\subseteq P$.  Given
$x\leq z$ in  $P$, we denote by $[x,z]\sub Q$ the subposet
$\{x,z\}\cup ([x,z]\cap Q)$ of $[x,z]$, and we write $[x,z)\sub Q$
for $[x,z]\sub Q\backslash\{z\}$.  We extend the definition 
of the M\"obius function of $Q$ to all intervals in $P$ by setting
\begin{equation}\label{eq:mobextend}
\mu\sub Q(x,z)=\mu\sub{[x,z]\sub Q}(x,z),
\end{equation} 
for all $x<z$ in $P$.
\begin{pro}\label{pro:mobius}
  Suppose that $P$ is a poset, $Q\subseteq P$ and $R=P\backslash Q$.
Then
\begin{equation}\label{eq:mobius}
\mu\sub P (x,z)\=\sum_{y\in [x,z)\sub Q}\mu\sub P(x,y)\mu\sub R(y,z),
\end{equation}
for all $x<z$ in $P$.
\end{pro}
\begin{proof}
  For each $y\in [x,z)\sub Q$, let $\cc_y$ denote the set of all
chains $C$ from $x$ to $z$ in $P$ such that $max{ (C\cap [x,z)\sub Q)}=y$.
From Theorem \ref{thm:hall}, we then have
$$
\mu\sub P(x,y)\mu\sub R(y,z)\= \sum_{C\in\cc_y}(-1)^{\ell (C)},
$$
for each such $y$. Since $\{\cc_y\co y\in [x,z)\sub Q\}$ is a partition
of the set of all chains from $x$ to $z$ in $P$, Equation \ref{eq:mobius}
thus follows from Theorem \ref{thm:hall}.
\end{proof}
\begin{pro}\label{pro:Pzero}
  Suppose that $P$ is a poset with minimum element $x$, and let
$\widehat P=\{y\in P\co\mu\sub P (x,y)\neq 0\}$.  
Then $\mu\sub{\widehat P}(x,z)
=\mu\sub P(x,z)$, for all $z\in\widehat P$.
\end{pro}
\begin{proof}
  The proof is immediate from the recursive definition of M\"obius function.
\end{proof}

\section{The free product of matroids}
We write $M(S)$ to indicate that $M$ is a matroid with underlying set
$S$. We denote by $\rho\sub M$ the rank function of $M(S)$, and write
$\rho (M)$ for the rank $\rho\sub M(S)$ of $M$.  The {\it free
product} of matroids $M(S)$ and $N(T)$ on disjoint sets $S$ and $T$
is the matroid $M\frp N$ on the union $S\cup T$ whose bases are those sets
$B\subseteq S\cup T$ of cardinality $\rho (M)+\rho (N)$
such that $B\cap S$ is independent in $M$ and $B\cap T$
spans $N$\footnote[2]{ This elegant characterization of free product
  in terms of bases is due to one of the referees of our paper
  \cite{crsc:uft}}.  The free product operation was introduced by the
current authors in \cite{crsc:fpm}, where it was used to prove the
conjecture of Welsh \cite{we:bnm} that $f_{n+m}\geq f_n\cdot f_m$,
where $f_n$ is the number of distinct isomorphism classes of matroids
on an $n$-element set.  In \cite{crsc:uft}, we studied the free
product in detail; in particular we showed that this operation, which
is noncommutative, is associative, and respects matroid duality in the
sense that
\begin{equation}
  \label{eq:dual}
  (M\frp N)\dual\= N\dual\frp M\dual,
\end{equation}
for all matroids $M$ and $N$.  We also characterized, in terms of
cyclic flats, those matroids which are irreducible with respect to
free product, and proved the following unique factorization theorem:
\begin{thm}\label{thm:uft}
    If $M_1\frp\cdots\frp M_k\isom N_1\frp\cdots\frp N_r$,
where each $M_i$ and $N_j$ is irreducible with respect to free product,
then $k=r$ and $M_i\isom N_i$, for $1\leq i\leq k$. 
\end{thm}
We gave
in \cite{crsc:uft} a number of cryptomorphic definitions of free
product, one of the most useful of which is the following proposition.
\begin{pro}\label{pro:rank}
The rank function of $M(S)\frp N(T)$ is given by
$$
\rho\sub{M\frp N}(A) 
\= \rkm{A\cap S} + \rkn{A\cap T} +\min\{\rkm S-\rkm{A\cap S},
\nln{A\cap T}\}, 
$$
for all $A\subseteq S\cup T$. 
\end{pro}
It is worth contrasting the above formula with that for the rank
function of the direct sum $M(S)\oplus N(T)$:
\begin{equation}
  \label{eq:sumrank}
  \rho\sub{M\oplus N}(A) 
\= \rkm{A\cap S} + \rkn{A\cap T},
\end{equation}
for all $A\subseteq S\cup T$. 

We refer to a matroid as {\it irreducible} if it is irreducible with
respect to free product, and {\it reducible}, otherwise.  
\begin{exa}\label{exa:freeprods}
We denote by $Z(a)$ and $I(a)$, respectively, the matroids consisting
of a single loop and single isthmus on $\{a\}$.  For
any set $S=\{a_1,\dots,a_n\}$, and $0\leq r\leq n$, 
the free product
$Z(a_1)\frp\cdots\frp Z(a_{n-r})\frp I(a_{n-r+1})\frp\cdots\frp I(a_n)$
is equal to $Z(a_1)\oplus\cdots\oplus Z(a_{n-r})\oplus I(a_{n-r+1})
\oplus\cdots\oplus I(a_n)$, the direct sum of $n-r$ loops and $r$ isthmi, 
while
$I(a_1)\frp\cdots\frp I(a_r)\frp Z(a_{r+1})\frp\cdots\frp Z(a_n)$
is the uniform matroid $U_{r,n}(S)$ of rank $r$ on $S$.  The matroid
$I(a)\frp Z(b)\frp I(c)\frp Z(d)$ is a three-point line, with
one point, $ab$, doubled.
\end{exa}

\begin{exa}
A matroid consisting of a single loop or isthmus is irreducible, 
and no matroid of size two or three is irreducible.  Up to isomorphism,
the unique irreducible matroid on four elements is the pair of 
double points $\un 12\oplus\un 12$, and on five elements the only
irreducibles are $\un 13\oplus\un 12$ and its dual $\un 23\oplus\un 12$.
\end{exa}

For any finite set $S$, we denote by $\cw\sub S$ the collection
of all matroids having $S$ as ground set.  The set $\cw\sub S$ is
partially ordered by the ({\it rank-preserving}) 
{\it weak order}, in which $M\leq N$ means that 
every basis for $M$ is also a basis for $N$ or, equivalently,
$M$ and $N$ have the same rank and the identity map on $S$ is
a weak map from $N$ to $M$.  The second inequality of the following 
result was Proposition 4.2 in \cite{crsc:uft}.
\begin{pro}\label{pro:univ}
For all matroids $M(S)$, and all $U\subseteq S$, the relation 
$$
M|U\oplus M/U\,\leq\, M\,\leq\, M|U\frp M/U
$$
holds in $\cw\sub S$.
\end{pro}
\begin{proof}
Let $V=S\backslash U$. The bases of $M|U\oplus M/U$ are those
subsets $B$ of $S$ such that $B\cap U$ is a basis for $M|U$ and
$B\cap V$ is a basis for $M/U$, which is the case if and only if
$B$ is a basis for $M$ such that $B\cap U$ spans $M|U$.  Hence
any basis for $M|U\oplus M/U$ is also a basis for $M$, and so
$M|U\oplus M/U\leq M$.
  
 By definition of free product, 
$B\subseteq S$ is a basis of $M|U\frp M/U$ if and only if
$|B|=\rho (M)$, the set $B\cap U$ is independent in $M$ and 
$B\cap V$ spans $M/U$. Now $\rho\sub{M/U}(B\cap V)=\rkm{B\cup U}-\rkm U$ 
and $\rho (M/U) = \rho (M) - \rkm U$, so $B\cap V$ spans $M/U$ if and only
$\rkm{B\cup U}=\rho (M)$, that is, if and only if $B\cup U$ spans $M$.  
Hence any basis for $M$ is also a basis for $M|U\frp M/U$, and so
$M\leq M|U\frp M/U$.
\end{proof}

The following result is one of the keys to understanding the
coproduct, and thus the primitive elements, of the matroid-minor
coalgebra.
\begin{pro}\label{pro:interval}
  For all matroids $M(S)$ and $N(T)$, with $S$ and $T$ disjoint, the set
$\{ L\,\co\, \text{$L|S = M$ and $L/S=N$}\}$
is equal to the interval $[M\oplus N, M\frp N]$ in $\cw\sub{S\cup T}$.
\end{pro}
\begin{proof}
If $L(S\cup T)$ is such that $L|S=M$ and $L/S=N$ then, by Proposition 
\ref{pro:univ}, $L\in [M\oplus N, M\frp N]$.
Conversely, suppose that $L\in [M\oplus N, M\frp N]$, so that
$\rho\sub{M\oplus N}(A)\leq\rho\sub L (A)\leq\rho\sub{M\frp N}(A)$,
for all $A\subseteq S\cup T$.
  If $A\subseteq S$ then, $\rkn{A\cap T}=\nln{A\cap T}=0$,
and so by Proposition \ref{pro:rank} and Equation \ref{eq:sumrank}, we have
$\rho\sub{M\oplus N}(A)=\rho\sub{M\frp N}(A)=\rkm A$.  Hence
$\rho\sub{L|S}(A)=\rkl A = \rkm A$, and therefore $L|S = M$.  

Now if $A\subseteq T$, then $\rho\sub{M\oplus N}(A\cup S) =
\rho\sub{M\frp N}(A\cup S) = \rkm S +\rkn A$, and thus
$\rkl{A\cup S}=\rkm S +\rkn A$.  Since
$\rho\sub{L/S}(A)= \rkl{A\cup S}-\rkl S$, and 
$\rkl S = \rkm S$, it follows that $\rho\sub{L/S}(A)=\rkn A$.
Hence $L/S=N$.
\end{proof}

\section{The matroid-minor coalgebra}
Let \cu\ be the set of all finite subsets of some fixed infinite
set, and let $K$ be a field.  Denote by \cw\ the set of all matroids
$M(S)$ whose ground set $S$ belongs to \cu, and denote by $\cw\sub
+$ the set of all nonempty $M\in\cw$.  We give \cw\ the (rank-preserving) 
weak order, in which $M\leq N$ means that $M$ and $N$ have the same ground
set $S$ and $M\leq N$ in $\cw\sub S$; hence \cw\ is the disjoint 
union of the posets $\cw\sub S$, for $S\in\cu$.  For sets $S$, $U$
and $V$, we write $U+V=S$ to indicate that $U\cup V = S$ and
$U\cap V=\emptyset$.

Let \frmod\cw\ denote the free $K$-vector space having basis \cw.  Since
$\{\cws\co S\in\cu\}$
is a partition of \cw, we have the direct sum decomposition
$$
\frmod\cw\= \bigoplus_{S\in\cu}\frmod\cws;
$$
hence the vector space \cw\ is graded by the set \cu, with
each homogeneous component $\frmod\cws$ finite-dimensional.

We define a pairing $\pr\cdot\cdot$ on $\frmod\cw$
by setting \pr MN\ equal to the Kronecker delta $\delta\sub{M,N}$, 
for all $M,N\in\cw$, and thus identify $\frmod\cw$ with the
graded dual space
$$
\frmod\cw\dual\= \bigoplus_{S\in\cu}\frmod\cws\dual.
$$

Let $C$ be the $K$-coalgebra on $\frmod\cw$
with coproduct $\delta$ and counit $\epsilon$
determined by
$$
\delta (M)\= \sum_{A\subseteq S}M|A\o M/A
\spandsp
\epsilon (M)\= 
\begin{cases}
1, & \text{if $S=\emptyset$},\\
0, & \text{otherwise,}
\end{cases}
$$
for all $M(S)\in\cw$.  Let $C\dual$ be the $K$-algebra on
\frmod\cw\ dual to $C$; the product of $C\dual$
is thus determined by 
$$
\pr{P\cdot Q}M\=\pr{\delta (M)}{P\o Q},
$$
for all $M,P,Q\in\cw$, and the unit element of $C\dual$ is
$\epsilon$.
 We remark that, even though the underlying
vector space of $C\dual$ is the graded dual of that of $C$, we do not
refer to $C\dual$ as the graded dual algebra of $C$, because $C$ is
not \cu-graded as a coalgebra.  In fact, since \cu\ has no
given monoid structure, the concept of \cu-graded coalgebra is
meaningless.

\begin{pro}\label{pro:prodform}
  In the algebra $C\dual$, the product of matroids $M$ and $N$
on disjoint ground sets is given by 
$$
M\cdot N\;\,\=\!\!\!\!\!\!\!
{\displaystyle\sum\sub{M\oplus N\leq L\leq M\frp N}}\!\!\!\!\!\!\!L. 
$$
If $M$ and $N$ are not disjoint, then $M\cdot N= 0$ in $C\dual$.
\end{pro}
\begin{proof}
First, observe that
  \begin{align*}
    M\cdot N & \= \sum_{L\in\cw}\,\pr{M\cdot N}L\, L\\
&\= \sum_{L\in\cw}\,\pr{M\o N}{\,\delta (L)}\, L.
  \end{align*}
If $M$ and $N$ are not disjoint, the latter sum is empty,
and hence $M\cdot N=0$.  If $M$ and $N$ are disjoint then,
by Proposition \ref{pro:interval}, 
$\pr{M\o N}{\delta (L)}$ is equal to one whenever
$M\oplus N\leq L\leq M\frp N$ in \cw, and is zero otherwise.
\end{proof}
More generally, the product of matroids $M_1,\ldots M_k$ in
$C\dual$ is equal to
\begin{equation}\label{eq:kprodform}
{\displaystyle\sum\sub{M_1\oplus\cdots\oplus M_k
\leq L\leq M_1\frp\cdots\frp M_k}}\!\!\!\!\!\!\!\!\!\!\!\!\!
\!\!\!\!\!\!\!L,
\end{equation}
if the set $\{M_1,\dots,M_k\}$ is pairwise disjoint and is zero
otherwise.
For all $S\in\cu$, we let $\pi\sub S\co\frmod\cw\rta\frmod\cws$
denote the natural
projection, and for any subset $X$ of $\frmod\cw$, we write $X\sub S$
for $\pi\sub S(X)$.  In particular, we write $C\sub S$ for
$\frmod\cws$ when viewed as a subspace of the coalgebra $C$, and
similarly for the algebra $C\dual$.
\begin{pro}\label{pro:grade}
  The coproduct $\delta$ of $C$ satisfies
$$
\delta (C\sub S)\subseteq\bigoplus\sub{U+V=S} C\sub U\o C\sub V,
$$
for all $S\in\cu$.
\end{pro}
\begin{proof}
  The result follows immediately from the definition
of $\delta$.
\end{proof}
We remark that
proposition \ref{pro:grade} may be stated alternatively
as
\begin{equation}\label{eq:grade}
\delta(\pi\sub S(x))\=\sum\sub{U+V=S}
(\pi\sub U\o\pi\sub V)\delta (x)
\end{equation}
for all $x\in C$ and $S\in\cu$.

Proposition \ref{pro:grade}, together with the fact that $\epsilon
(C\sub S)=0$, for all $S\neq\emptyset$ says that $C$ is something very
much like a \cu-graded coalgebra; the problem, as we indicated above,
is that the disjoint union operation $+$ is only partially defined on
\cu\ and so does not equip \cu\ with a monoid structure.  To make
precise the sense in which $C$ is a ``generalized graded'' coalgebra,
consider the {\it partial monoid algebra} \frmod\cu, with product
determined by
$$
ST\=
\begin{cases}
  S\cup T & \text{if $S\cap T = \emptyset$,}\\
0 & \text{otherwise,}
\end{cases}
$$
and coproduct $\delta (S)=S\o S$, for all $S,T\in\cu$.  Then
\frmod\cu\ is a Hopf algebra, and the map $\psi\co C\rta\frmod\cu\o
C$, determined by $x\mapsto S\o x$, for all $x\in C\sub S$, is a
coaction, making $C$ a $\frmod\cu$-comodule coalgebra.  In the
following situation: $\cu$ a monoid, $C\sub S$ a subspace of a
coalgebra $C$, for all $S\in\cu$ and $\frmod\cu$ the monoid (Hopf)
algebra of \cu, the above map $\psi$ being well-defined and making $C$ a
$\frmod\cu$-comodule coalgebra is equivalent to $C$ being \cu-graded
as a coalgebra, with homogeneous components $C\sub S$, for all
$S\in\cu$.

\section{Primitive elements in the matroid-minor coalgebra}\label{sec:prim}

We denote by
\ci\ and \ccr, respectively, the collections of all irreducible and
all reducible matroids belonging to \cw.  For any $M\in\cw$, we denote
by $\cw\sub M$ the order filter $\{ N\in\cw\co N\geq M\}$ of \cw\ and
define the following subposets of \cwm:
$$
\cim \= (\ci\cap\cwm)\cup\{M\}\spandsp
\crm \=(\ccr\cap\cwm)\cup\{M\}.
$$
For any $S\in\cu$, we set
$$\cis= \ci\cap\cws\spandsp\crs=\ccr\cap\cws.
$$

Let $C\sub + =\ker\epsilon$, and let $\bar\delta\co C\sub +\rta C\sub
+\o C\sub +$ be the map determined by $\delta (x)\= 1\o x + x\o 1 +
\bar\delta (x)$, for all $x\in C\sub +$, where $1$ denotes the empty
matroid.  Then $C\sub +$ has basis $\cw\sub +$, and
$\bar\delta (x)$ satisfies
$$
\bar\delta (M)\= \sum_{\stackrel{A\subseteq S}
{\scriptscriptstyle A\neq\emptyset, S}} M|A \o M/A,
$$
for all $M\in\cw\sub +$.  We write $P(C)$ for the subspace of
primitive elements of $C$; hence $P(C) = \{ x\in C\sub +\co\,
\bar\delta (x)=0\}$.
\begin{pro}\label{pro:primgrade}
  The space of primitive elements $P(C)$ 
respects the grading of $C$ by $\cu$,
that is,
$$
P(C)\sub S\= P(C)\cap C\sub S,
$$
for all $S\in\cu$.
\end{pro}
\begin{proof}
  For any $S\in\cu$ we have $P(C)\cap C\sub S =\pi\sub S(P(C)\cap
  C\sub S) \subseteq P(C)\sub S$.  On the other hand, if $x\in
  P(C)\sub S$, then $x=\pi\sub S(y)$ for some $y\in P(C)$ and thus,
  by Equation \ref{eq:grade},
\begin{align*}
  \delta (x)& \= \sum\sub{U+V=S}(\pi\sub U\o\pi\sub V)\delta (y)\\
&\= \sum\sub{U+V=S}\left(\pi\sub U(y)\o \pi\sub V(1) + 
\pi\sub U(1)\o\pi\sub V(y)\right).
\end{align*}
Since
$$
\pi\sub V(1) =
\begin{cases}
  1 & \text{if $V = \emptyset$,}\\
0 & \text{otherwise,}
\end{cases}
$$
we thus have
\begin{align*}
  \delta(x)&\=\pi\sub S(y)\o 1 + 1\o\pi\sub S(y)\\
&\= x\o 1 + 1\o x.
\end{align*}
Hence $P(C)\sub S\subseteq P(C)\cap C\sub S$.  
\end{proof}

An alternative way of stating Proposition \ref{pro:primgrade}
is the following:
\begin{equation}\label{eq:primdecomp}
P(C)\= \bigoplus_{S\in\cu}P(C)\cap C\sub S.
\end{equation}

Denote by $C\dual\sub +$ the ideal in $C\dual$ consisting of all
elements $x$ such that $\langle x\,, 1\rangle = 0$.  Note that
$C\dual\sub+$ has basis $\cw\sub +$, and that the ideal $(C\dual\sub
+)^2$ of $C\dual$ is spanned by the set of all products $P\cdot Q$
such that $P,Q\in\cw\sub +$.  For any subset $X$ of $\frmod\cw$,
we define $X^\perp =\{y\in\frmod\cw\co\text{$\pr yx=0$, for all
$x\in X$}\}$. The following proposition is a standard
result about connected coalgebras.
\begin{pro}\label{pro:prim1}
The subspace of primitive elements of $C$ is given by
$P(C) = C\sub +\cap [(C\dual\sub +)^2]^\perp$.
\end{pro}
\begin{proof}
An element $x$ of $C$ belongs to $P(C)$ if and only if
$\epsilon (x)=0$ and $\bar\delta (x)=0$, that is, if and only
if $x\in C\sub+$ and $\pr{P\o Q}{\bar\delta (x)}=0$, for
all $P,Q\in\cw\sub +$.  But, for nonempty $P$ and $Q$, we have
\begin{align*}
\pr{P\o Q}{\,\bar\delta (x)} & \= \pr{P\o Q}{\,\delta (x)}\\
&\=\pr{P\cdot Q}{\,x}.
\end{align*}
Hence, the condition that  $\pr{P\o Q}{\bar\delta (x)}=0$, for
all $P,Q\in\cw\sub +$, is equivalent to $x$ belonging to 
$[(C\dual\sub +)^2]^\perp$.
\end{proof}

\begin{cor}\label{cor:primgrade}
  For all nonempty $S\in\cu$ we have
$$
P(C)\sub S\= ([(C\dual\sub +)^2]\sub S)^\perp,
$$
while $P(C)_{\emptyset}=\{0\}$.
\end{cor}
\begin{proof}
  The result follows directly from Proposition \ref{pro:prim1}
and the direct-sum decomposition \eqref{eq:primdecomp}.
\end{proof}
Untangling the notation of Corollary \ref{cor:primgrade}, we see
that $P(C)\sub S$ consists of all $x\in C\sub S$ such that
$\pr x{P\cdot Q}=0$, for all matroids $P(U)$ and $Q(V)$
with $U$ and $V$ nonempty and $U+V=S$.

The following result streamlines the brute-force method of
finding primitive elements in $C$.
\begin{pro}\label{pro:knead}
  The ideal $(C\dual\sub +)^2$ of $C\dual$ is spanned by the set \cp\
of all products $P\cdot Q$ such that $P$ is irreducible
(with respect to free product) and $Q$ is nonempty.
\end{pro}
\begin{proof}
We prove the result by using (weak order) induction on $P$
to show that if $P\cdot Q\in (C\dual\sub +)^2$, then 
$P\cdot Q\in\frmod\cp$.

For the base case, suppose that $P\cdot Q\in (C\dual\sub +)^2$ and
that $P=P(S)$ is minimal in the weak order. Then $P$ is the direct
sum of $k$ loops and $r$ isthmi, where $k$ is the nullity and $r$
is the rank of $P$, and so, in the notation of Example \ref{exa:freeprods},
we have
\begin{align*}
P(S)& \= 
Z(a_1)\oplus\cdots\oplus Z(a_k)\oplus I(a_{k+1})\oplus\cdots\frp I(a_{k+r})\\
&\=
Z(a_1)\frp\cdots\frp Z(a_k)\frp I(a_{k+1})\frp\cdots\frp I(a_{k+r}),
\end{align*}
which, by \eqref{eq:kprodform}, implies that $P$ is equal to 
the product
$$
Z(a_1)\cdots Z(a_k)\cdot I(a_{k+1})\cdots I(a_{k+r}),
$$
in $C\dual$. Hence $P\cdot Q\in\frmod\cp$.  

Now suppose that $P\cdot Q\in (C\dual\sub +)^2$, and that 
$L\cdot Q\in\frmod\cp$,
for all $L<P$.  If $P$ is irreducible, we're done; otherwise, write
$P$ as $M\frp N$, where $M$ is irreducible.  By Proposition 
\ref{pro:prodform}, we have
$$
M\cdot N\cdot Q\,\,\,\,=\!\sum\sub{M\oplus N\leq L\leq P}\!\!\! L\cdot Q,
$$
and hence
$$
P\cdot Q\,\=\, M\cdot N\cdot Q\,\,\,\, -\!
\sum\sub{M\oplus N\leq L< P} \!\!\! L\cdot Q,
$$
which, by induction, belongs to $\frmod\cp$.
\end{proof}
Of course, the ``right-handed'' version of Propositon \ref{pro:knead} 
also holds: the ideal $(C\dual\sub +)^2$ of $C\dual$ is spanned by the 
set of all products $Q\cdot P$ such that $P$ is irreducible
and $Q$ is nonempty.
\begin{pro}\label{pro:primdim}
The inequality $\dim P(C)\sub S\leq |\cis|$ holds
for all $S\in\cu$. 
\end{pro}
\begin{proof}
Define a map $\alpha\co\frmod{\crs}\rta (C\dual\sub +)^2$ as follows:
for each reducible $M(S)$, choose a sequence of irreducible
matroids $M_1,\dots,M_k$ such that $M=M_1\frp\cdots\frp M_k$
(recall that the sequence $M_1,\dots, M_k$ is uniquely determined
by $M$ only up to isomorphism), then set $\alpha (M)$ equal to
the product $M_1\cdots M_k$ in $C\dual$.   Clearly $\text{\rm im}\,
\alpha\subseteq [(C\dual\sub +)^2]\sub S$.
By Proposition \ref{pro:prodform},
we know that $\alpha (M)$ is equal to $M$ plus matroids that
are less than $M$ in the weak order.  It follows that $\alpha$
is injective, and so $|\crs|\leq\dim [(C\dual\sub +)^2]\sub S$.  Hence,
by Corollary \ref{cor:primgrade}, we have
\begin{align*}
\dim P(C)\sub S& \= \dim ([(C\dual\sub +)^2]\sub S)^\perp\\
&\= |\cws|-\dim [(C\dual\sub +)^2]\sub S\\
&\,\leq |\cws| -|\crs|\\
&\= |\cis|.
\end{align*}
\end{proof}
We will see shortly, in Theorem \ref{thm:primbasis1}, that
the inequality in the above proposition is in fact
an equality.   Recall that a {\it free separator}
of a matroid $M(S)$ is a subset $U$ of $S$ such that $M$ factors
as $M=M|U \frp M/U$.
\begin{pro}
  For all $S\in\cu$ and $U\subseteq S$, the
map $\varphi\sub U\co\cw\sub S\rta\cw\sub S$ given by
$$
\varphi\sub U (M)\= M|U \frp M/U,
$$
for all $M\in\cw\sub S$, is a closure operator.
A matroid $M$ is $\varphi\sub U$-closed if and only if
$U$ is a free separator of $M$.  
\end{pro}
\begin{proof}
Let $M$ be a matroid on $S$, and $U\subseteq S$. 
It is immediate from the definitions that
$U$ is a free separator of $M$ if and only if
$M=\varphi\sub U (M)$.  To show that $\varphi\sub U$ is
a closure operator, we first note that, since $U$ is a free separator of 
$\varphi\sub U(M)$, we have $\varphi\sub U (\varphi\sub U(M)) =
\varphi\sub U(M)$; also, the inequality $M\leq\varphi\sub U (M)$ follows
from Proposition \ref{pro:univ}.
It remains to show that, if $M\leq N$ in
$\cw\sub S$, then $\varphi\sub U (M)\leq\varphi\sub U (N)$.  
As we noted in the proof of Proposition \ref{pro:univ},
The bases of $M|U\frp M/U$ are those subsets $B$ of $S$
with $|B|=\rho (M)$ such that $B\cap U$ is independent in $M$
and $B\cup U$ spans $M$.  Now if $N\geq M$ and $B\cap U$ is
independent in $M$, then it is also independent in $N$; since
$\rho (N)=\rho (M)$, it thus follows that any basis for $M|U\frp M/U$ is
also a basis for $N|U\frp N/U$, that is, $\varphi\sub U (M)\leq\varphi\sub U
(N)$.
\end{proof}

\begin{lem}\label{lem:clch}
  For all matroids $P(U)$, $Q(V)$ and $M(S)$, with $S=U+V$,
we have $\varphi\sub U (M) = P\frp Q$ if and only if $P\oplus Q\leq
M\leq P\frp Q$ in \cws.
\end{lem}
\begin{proof}
Since $(P\frp Q)|U=P$ and $(P\frp Q)/U=Q$, for all matroids $P(U)$ and 
$Q(V)$, it follows that $P(U)\frp Q(V) = \varphi\sub U (M) = M|U\frp M/U$ 
if and only if $M|U = P$ and $M/U = Q$. By Proposition \ref{pro:interval}
this is the case if and only if $P\oplus Q\leq M\leq P\frp Q$.
\end{proof}

We define bases $\{w\sub M\co M\in\cw\}$ and $\{r\sub M\co M\in\cw\}$
for the matroid-minor coalgebra $C$ by setting
$$
w\sub M \= \sum_{N\in\cwm}\mu\sub\cw (M,N)\,N
\spandsp
r\sub M \= \sum_{N\in\crm}\mu\sub\ccr (M,N)\,N,
$$
for all $M\in\cw$.  We have written here $\mu\sub\ccr$ for the
M\"obius function of the poset $\crm$; equivalently, we are
using the convention \eqref{eq:mobextend} to extend the definition
of the M\"obius function of \ccr\ to arbitrary intervals in \cw.
Note that $\{w\sub M\co M\in\cw\}$ and $\{r\sub
M\co M\in\cw\}$ are indeed bases for $C$ since, by M\"obius inversion,
we have
$$
M \= \sum_{N\in\cwm}\!\!w\sub N
\= \sum_{N\in\crm}\!\!r\sub N,
$$
for all $M\in\cw$.  The following proposition shows us how
to change between the bases  $\{w\sub M\co M\in\cw\}$ and $\{r\sub
M\co M\in\cw\}$.
\begin{pro}
  For all $M\in\cw$
$$
w\sub M \= \sum_{N\in\cim}\mu\sub\cw (M,N)\,r\sub N
\spandsp
r\sub M \= w\sub M - \!\!\sum_{N\in{\cim\backslash M}}
\mu\sub\ccr (M,N)\,w\sub N,
$$
\end{pro}
\begin{proof}
By definition of $r\sub N$, $w\sub M$, and Proposition
\ref{pro:mobius}, we have
\begin{align*}
    \sum_{N\in\cim}\mu\sub\cw (M,N)\,r\sub N & \=
\sum_{N\in\cim}\sum_{Q\in\ccr\sub N}
\mu\sub\cw (M,N)\mu\sub\ccr (N,Q)\,Q\\
&\= M \+ \sum_{Q\in\cim\backslash M}\mu\sub\cw (M,Q)\\
& \+
\sum_{Q\in\crm\backslash M}\sum_{N\in [M,Q)\sub\ci}
\mu\sub\cw (M,N)\mu\sub\ccr (N,Q)\, Q\\
& \= \sum_{Q\in\cwm}\mu\sub\cw (M,Q)\,Q\\
&\= w\sub M, 
\end{align*}
thus establishing the first equality.  For the second, we compute
\begin{align*}
  r\sub M & \= \sum_{N\in\crm}\mu\sub\ccr (M,N)\, N\\
& \= \sum_{N\in\crm}\sum_{Q\in\cw\sub N}\mu\sub\ccr (M,N)\,w\sub Q\\
&\= w\sub M \+ \sum_{Q\in\crm\backslash M}\sum_{N\in [M,Q]\sub R}
\mu\sub\ccr (M,N)\,w\sub Q\\
&\+ \sum_{Q\in\cim\backslash M}\sum_{N\in [M,Q)\sub R}
\mu\sub\ccr (M,N)\,w\sub Q.\\
\end{align*}
By the recursive definition of M\"obius function it follows that
$$
\sum_{N\in [M,Q]\sub\ccr}\mu\sub\ccr (M,N)\= 0\spandsp
\sum_{N\in [M,Q)\sub\ccr}\mu\sub\ccr (M,N)\= -\mu\sub\ccr (M,Q),
$$
for all $Q>M$ in \cw; thus the result follows.
\end{proof}

We now come to the first of our main results.
\begin{thm}\label{thm:primbasis1}
The set $\{w\sub M\co M\in\ci\}$ is a homogeneous basis for $P(C)$.
\end{thm}
\begin{proof}
It is clear from the definition that the $w\sub M$
are homogeneous with respect to the \cu-grading of $C$.  
Now, if $M$ is irreducible, then it is 
nonempty and thus $\epsilon (w\sub M) = 0$. 
Suppose that $M(S)$ is irreducible and that $P(U)$ and $Q(V)$ are
nonempty matroids.  By Proposition \ref{pro:prodform}, we have
$$
\pr{P\cdot Q}{w\sub M}\=\!\!\!\!
\sum_{\stackrel{N\in\cwm}{\scriptscriptstyle
P\oplus Q\leq N\leq P\frp Q}}\!\!\!\!\!\mu\sub\cw (M,N).
$$
If it is not the case that
$M\leq P\frp Q$,
then the sum is empty; otherwise, according to Lemma \ref{lem:clch},
 it is given by
$$
\sum_{\stackrel{N\in\cwm}{\scriptscriptstyle\varphi\sub U(N)=P\frp Q}
}\mu\sub\cw (M,N).
$$
Since $U$ is a nonempty proper subset of $S$ and $M$ is
irreducible, $U$ is not a free separator of $M$, and so $M$ is not
$\varphi\sub U$-closed.  Thus, by Theorem \ref{thm:closure},
the above sum is zero, and so it follows from Proposition
\ref{pro:prim1} that $w\sub M$ is primitive in $C$.

Since $\{w\sub M\co M\in\cw\}$ is a basis for $C$, the set
$\{w\sub M\co M\in\cis\}$ is linearly independent in $P(C)\sub S$ and
thus, by Proposition \ref{pro:primdim}, is a
basis for $P(C)\sub S$. It follows from the direct-sum decomposition 
\eqref{eq:primdecomp} that $\{w\sub M\co M\in\ci\}$ is a basis for $P(C)$.
\end{proof}

\begin{thm}\label{thm:primbasis2}
The set $\{r\sub M\co M\in\ci\}$ is a homogeneous
basis for $P(C)$.  Furthermore, if an element $x$ of $P(C)$ has
the form $M + y$,
where $M$ is an irreducible matroid and $y$ is a linear
combination of reducible matroids, the $x=r\sub M$.
\end{thm}
\begin{proof}
Observe that, for any $M(S)$, and $U\subseteq S$, the closure
operator $\varphi\sub U$ on \cwm\ satisfies $\varphi\sub U(\crm)\subseteq
\crm$, and thus restricts to a closure operator on \crm.  The proof
that $\{r\sub M\co M\in\ci\}$ is a basis for $P(C)$ thus parallels
that of Theorem \ref{thm:primbasis1}, with
the poset \crm\ used in place of \cwm.
Uniqueness follows immediately from the fact
that the set \ci\ of irreducible matroids is linearly
independent in $C$.
\end{proof}

\begin{exa}\label{exa:D}
  Consider the irreducible matroid $D=U_{1,2}(a,b)\oplus U_{1,2}(c,d)$,
consisting of two double points $ab$ and $cd$.  The poset $\cw\sub D$
consists of $D$, the four-point line $Q=U_{2,4}(a,b,c,d)$, and 
the matroids $P_1= I(a)\frp Z(b)\frp I(c)\frp Z(d)$ and $P_2= 
I(c)\frp Z(d)\frp I(a)\frp Z(b)$ (see Example \ref{exa:freeprods}), with
$D\leq P_1,P_2\leq Q$, and $P_1$ and $P_2$ are incomparable.
Since $\cw\sub D$ contains no irreducible matroids other than $D$, we have
$\cw\sub D=\ccr\sub D$; hence, by Theorem \ref{thm:primbasis1},
$$
w\sub D\= r\sub D\= D \,-\, P_1\, -\, P_2\+ Q
$$
is primitive in $C$.
\end{exa}

\begin{figure}\label{fig:lattice5}
  \includegraphics[scale=.8,bb=30 20 480 360]
{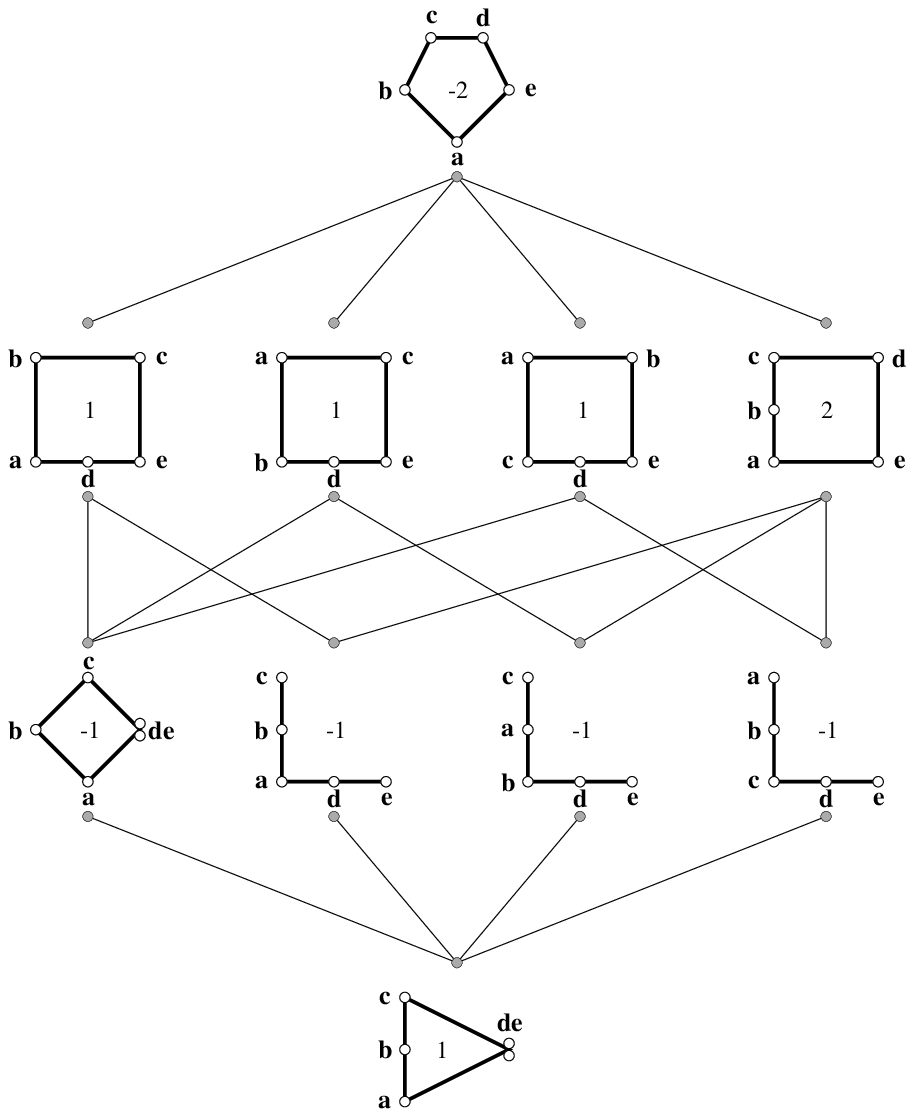}
  \caption{The poset $\cwm = \crm$, for 
$M=U_{2,3}(a,b,c)\oplus U_{1,2}(d,e)$}
\label{fig:littlelattice} 
\end{figure}

\begin{exa}\label{exa:lattice5}
The Hasse diagram of the poset $\cwm$, where $M$ is the direct sum of the 
three-point line $abc$ and the double point $de$, is shown in
Figure \ref{fig:littlelattice}. Since \cwm\ contains no irreducible matroids
other than $M$, we have $\cwm = \crm$. Each matroid $N$ shown in the 
diagram is labeled by the M\"obius function value 
$\mu\sub\cw (M,N)$, so the primitive element
$$
\includegraphics[scale=.7,viewport=49 50 475 100]
{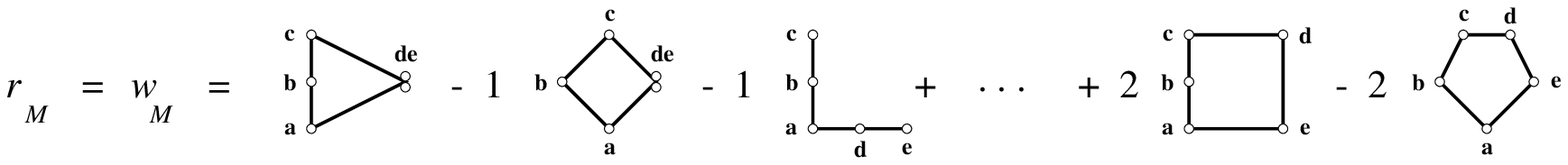}$$\\
can be read from the diagram.
\end{exa}

\begin{figure}\label{fig:prim7}
  \includegraphics[scale=.8,bb=80 20 480 720]
{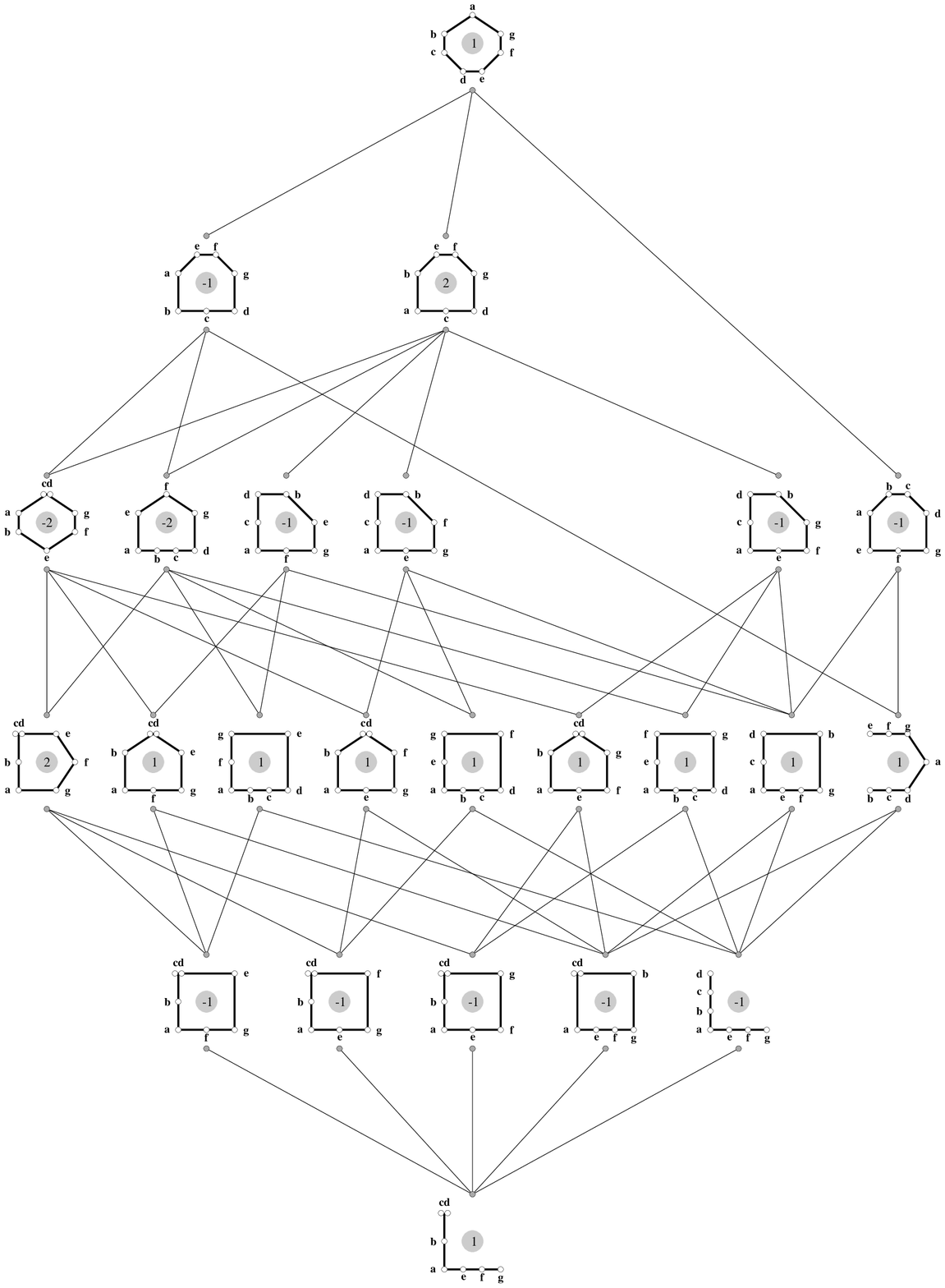}
  \caption{The poset ${\widehat\ccr}\sub M$,
for $M$ the matroid of Example \ref{exa:lattice7}}
  \label{fig:biglattice}
\end{figure}

\begin{exa}\label{exa:lattice7}
 Suppose that $M$ is the matroid shown below.\\
  \includegraphics[scale=.65,bb=-230 15 300 100]{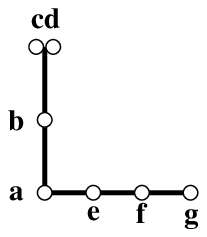}\\
The Hasse diagram of the poset ${\widehat\ccr}\sub M=\{ N\in\crm\co
\mu\sub\ccr (M,N)\neq 0\}$ is shown in Figure \ref{fig:biglattice}.
The matroids belonging to the set $\crm\backslash{\widehat\ccr}\sub M=
\{N\in\crm\co\mu\sub\ccr (M,N)=0\}$ are shown below.\\
 \includegraphics[scale=.7,bb=-27 15 350 200]{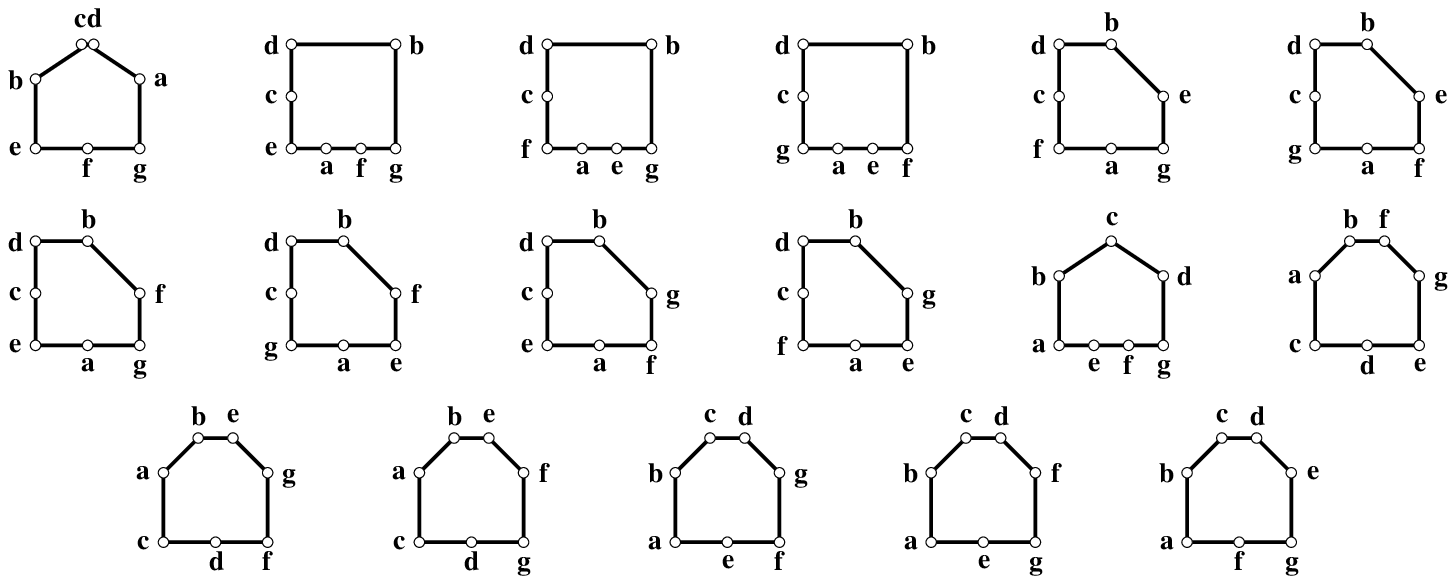}\\

Each matroid $N$ in the Hasse diagram is labeled by
the M\"obius function value $\mu\sub{\widehat\ccr}(M,N)$, 
which is equal to
$\mu\sub\ccr (M,N)$, by Proposition \ref{pro:Pzero}.    
Hence the primitive element
$$
\includegraphics[scale=.75,bb=20 45 410 110]
{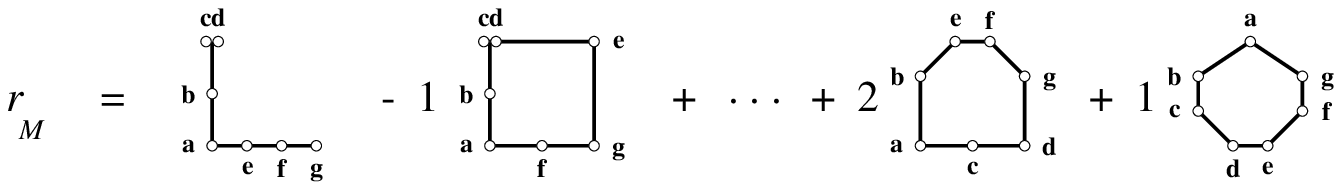}$$\\
can be read from the diagram. In contrast to the situation in
Examples \ref{exa:D} and \ref{exa:lattice5},
we have $\cwm\neq\crm$ here.  The set $\cwm\backslash\crm$
consists of $15$ matroids, namely, the following four
$$
\includegraphics[scale=.75,bb=0 240 390 310]
{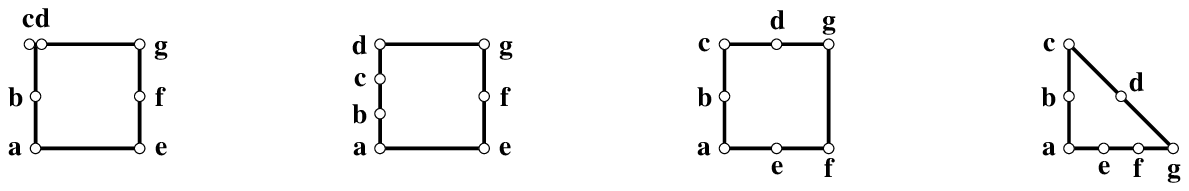}\\
$$
together with $11$ others, one of which is isomorphic to the
second of these, and five isomorphic to each of the
last two.
\end{exa}

\section{The matroid-minor Hopf algebra}

We now turn our attention to the matroid-minor Hopf algebra, which has
distinguished basis consisting of isomorphism classes of matroids,
rather than labeled matroids as we have been considering up to now.
We denote by \cm\ and \cmi, respectively, the set of all isomorphism
classes of matroids and the set of all isomorphism classes of
irreducible matroids.  (Recall that {\it irreducible} here means
irreducible with respect to free product.)  We write $\cl M$ for the
isomorphism class of a matroid $M$ by $\cl M$, and note that $\cm
=\{\cl M\co M\in\cw\}$.  The (rank-preserving) {\it weak order} on
$\cm$ is the coarsest order relation such that the map $\cw\rta\cm$,
taking a matroid to its isomorphism class, is order preserving; thus
$\cl M\leq\cl N$ in \cm\ if and only if there exist representatives
$M'\in\cl M$ and $N'\in\cl N$ such that $M'\leq N'$ in \cw.

The {\it matroid-minor Hopf algebra} $H$, first defined in
\cite{sc:iha}, is the free vector space $\frmod\cm$, with product
induced by direct sum, that is, $\cl M\cl N = \cl{M\oplus N}$, for all
matroids $M$ and $N$, and with
coproduct $\delta\sub H$ and counit
$\epsilon$ determined by
$$
\delta\sub H (\cl M)\= \sum_{A\subseteq S}\cl{M|A}\o\cl{M/A}
\spandsp
\epsilon (\cl M)\= 
\begin{cases}
1, & \text{if $S=\emptyset$},\\
0, & \text{otherwise,}
\end{cases}
$$
for all $M$.
Note that, as an algebra, $H$ is free commutative, generated by the
set $\cm_c$ of isomorphism classes of connected matroids, that is,
$H$ is isomorphic to the polynomial algebra $K[\cm_c]$.

As a means of determining
the coalgebra structure of $H$, we first introduce
a coalgebra, $L$, whose underlying vector space is also
\frmod\cm, but with coproduct dual to the product on \frmod\cm\ 
induced by the free product operation; that is, for all matroids $M$,
$$\delta\sub L (\cl M)= \sum\cl P\o\cl Q,
$$
where the sum is over all pairs $(\cl P, \cl Q)\in\cm\times\cm$ 
such that $P\frp Q\isom M$. The counit of $L$ is equal to that of $H$.
Associativity of
free product implies that $\delta\sub L$ is coassociative.
By Theorem \ref{thm:uft}, for every $\cl M\in\cm$, there
is a unique sequence $(\cl{M_1},\dots,\cl{M_k})$, for some $k\geq 0$,
with each $M_i$ irreducible, such that 
$M\isom M_1\frp\cdots\frp M_k$.  Denoting, for the moment,
each element of \cm\ by its corresponding sequence, or {\it word}, 
we see that the coalgebra $L$ has as
basis all words on the set \cmi\ of isomorphism classes of
irreducible matroids, and its coproduct can be written
\begin{equation}\label{eq:cofreecop}
\delta\sub L (\cl{M_1},\dots,\cl{M_k}) = 
\sum_{i=0}^k (\cl{M_1},\dots,\cl{M_i})\o
(\cl{M_{i+1}},\dots,\cl{M_k}).
\end{equation}
Hence $L$ is cofree, and $\cm\sub I$ 
is a basis for the space of primitive elements $P(L)$.

Our next result is a theorem from \cite{crsc:uft} which shows that $H$
and $L$ are isomorphic coalgebras, if the field $K$ has characteristic
zero.  In the proof we give here, which is dual to the one given in
that article, we use the following notation and terminology: If $M(S)$
is a matroid, and $A\subseteq B\subseteq S$, we denote by $M(A,B)$ the
minor $(M|B)/A=(M/A)|(B\backslash A)$ of $M$ determined by the
interval $[A,B]$ in the Boolean algebra of subsets of $S$.  For any
set $S$, an $S$-{\it chain} is a sequence of sets $(S_0,\dots , S_k)$
such that $S_{i-1}$ is strictly contained in $S_i$, for $1\leq i\leq
k$.  Given $S$-chains $C=(S_0,\dots , S_k)$ and $D=(T_0,\dots,T_\ell)$
such that $S_k=T_0$, we write $CD$ for the $S$-chain $(S_0,\ldots,
S_k=T_0,\ldots,T_\ell)$.  If $S$ and $T$ are disjoint sets,
$C=(S_0,\dots , S_k)$ an $S$-chain and $D=(T_0,\dots,T_\ell)$ a
$T$-chain with $T_0=\emptyset$, we denote by $C\cdot D$ the $(S\cup
T)$-chain $(S_0,\ldots, S_k=T_0\cup S_k,\ldots, T_\ell\cup S_k)$.
Given a matroid $M(S)$ and $S$-chain $C=(S_0,\dots,S_k)$, we write
$M(C)$ for the free product $M(S_0,S_1)\frp\cdots\frp M(S_{k-1},S_k)$,
We refer to an $S$-chain $C$ as $M$-{\it irreducible} if $S_0 =
\emptyset$, $S_k=S$ and each of the minors $M(S_{i-1},S_i)$ is
irreducible with respect to free product.  We write \ic M\ for the set
of all $M$-irreducible $S$-chains.

\begin{thm}[\cite{crsc:uft}]\label{thm:isom}
If the field $K$ has characteristic zero, then the 
map $\varphi\co H\rta L$, determined by
$$\varphi (\cl M) =
\sum\sub{C\in\ic M}\cl{M(C)},
$$
for all $M$, is a coalgebra isomorphism.
\end{thm}
\begin{proof}
We compute, for $M=M(S)$,
\begin{align*}
\delta\sub L (\varphi (\cl M))
&\= \sum\sub{C\in\ic M}\delta\sub L (\cl{M(C)})\\    
&\= \sum\sub{C\in\ic M}\sum\sub{C_1C_2=C} \cl{M(C_1)}\o\cl{M(C_2)},
\end{align*}
and
\begin{align*}
(\varphi\o\varphi)\delta\sub H (\cl M)
&\= \sum\sub{A\subseteq S}\varphi (M|A)\o\varphi (M/A)\\
&\= \sum\sub{A\subseteq S}
\sum\sub{\stackrel{D\in\ic{M|A}}{E\in\ic{M/A}}}
\cl{M(D)}\o\cl{M(E)}.
\end{align*}
It is readily verified that the map
$$
\bigcup_{A\subseteq S}\ic{M|A}\times\ic{M/A}\rta
\{ (C,A)\co\text{$C\in\ic M$ and $A\in C$}\},
$$
given by $(C_1,C_2)\mapsto (C_1\cdot C_2, A)$, for all $A\subseteq S$ and
$(C_1,C_2)\in\ic{M|A}\times\ic{M/A}$, is a bijection.  It thus
follows from the above computation that 
$\delta\sub L (\varphi (\cl M)) = (\varphi\o\varphi)\delta\sub H
(\cl M)$, and so $\varphi$ is a coalgebra map.

Proposition \ref{pro:univ} implies that $\cl M\leq \cl{M(C)}$ in \cm, 
for all $C\in\ic M$; since $K$ has characteristic zero, it thus follows that
$\varphi$ is a bijective and thus an isomorphism. 
\end{proof}

For all
$M\in\cm$, we write $p\sub M$ for $\varphi^{-1}(\cl M)\in H$. 
\begin{thm}[\cite{crsc:uft}]\label{thm:cofree}
If the field $K$ has characteristic zero, then the 
matroid-minor Hopf algebra $H$ is cofree, with $\{p\sub M\co \cl M\in\cmi\}$
a basis for the space $P(H)$ of primitive elements of $H$.
\end{thm}
\begin{proof}
  The result follows immediately from Theorem \ref{thm:isom} and
the preceeding discussion of the coalgebra $L$.
\end{proof}
Observe that, since $\cl M\leq\cl{M(C)}$,
for all matroids $M$ and chains $C\in\ic M$, we may write
$$
\varphi (\cl M)=\sum_{\cl N\geq\cl M}c(\cl M, \cl N)\, \cl N,
$$
where $c(\cl M, \cl N)$ is the cardinality of the
set $\{ C\in\ic M\co M(C)\isom N\}$.  Hence we may compute
$p\sub M = \varphi^{-1}(\cl M)$ by computing the inverse
the matrix whose entries are the numbers $c(\cl P,\cl Q)$,
for all $P\leq Q$ in \cwm.

\begin{exa}\label{exa:prim5uft}
Suppose that $M$ is the $5$-element matroid of Example
\ref{exa:lattice5}.  
The isomorphism classes of matroids
in \cwm\ are the following
$$
\includegraphics[scale=.65,bb=24 55 450 110]
{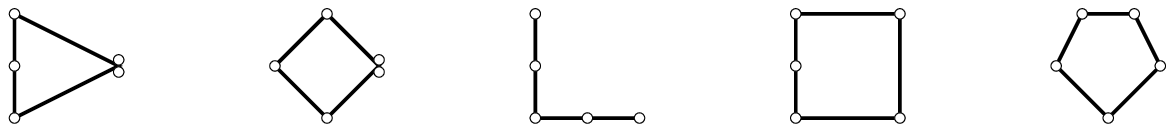}
$$\\
We refer to these isomorphism classes, in the order shown here,
as $e_1,\dots,e_5$.  Note that we have the following factorizations
into irreducibles:
$e_1=\cl M$ is irreducible, $e_2=IZIIZ$,
$e_3=ID$, 
$e_4=IIZIZ$ and $e_5=U_{3,5}=IIIZZ$,
where we have written $I$ and $Z$ for the isomorphism
classes of the single point and loop, and $D$ for the irreducible isomorphism
class $U_{1,2}\oplus U_{1,2}$
consisting of two double-points (see Example \ref{exa:D}), and we
have suppressed the symbol $\frp$ in writing free products
of isomorphism classes.  Denoting by $\Phi$ the matrix whose
entry in position $(i,j)$ is
$c(e_i,e_j)$, for $1\leq i\leq j$, we have

$$
\Phi\= \left(\begin{array}{crrrr}
1& 12 & 3 & 36 & 72\\
0 &12 & 0 & 24 & 84 \\
0 &0 & 1& 24& 96\\
0 &0 & 0 & 12 & 108 \\
0 &0 & 0 & 0 & 120\\
\end{array}
\right) 
\spandsp
\Phi^{-1}\=\left(\begin{array}{crrrr}
1& -1 & -3 & 5 & -2\\
0 & \frac{1}{12} & 0 & -\frac 16 & \frac{11}{120} \\
0 &0 & 1& -2& 1\\
0 &0 & 0 & \frac 1{12} & -\frac 3{40} \\
0 &0 & 0 & 0 & \frac 1{120}\\
\end{array}
\right).
$$\\
It follows that the basis element $p\sub M=\varphi^{-1}(\cl M)$
of $P(H)$ is given by
\begin{equation}\label{eq:prim5}
\includegraphics[scale=.65,bb=24 62 450 110]
{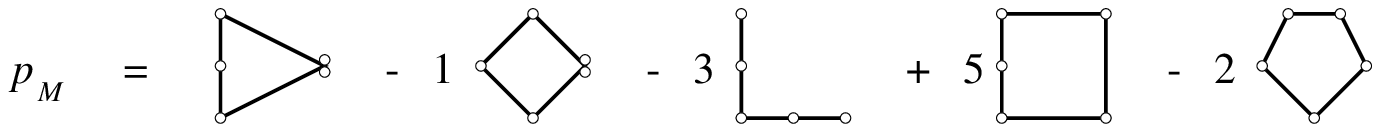}
\end{equation}\\
\end{exa}

The following result, which is analogous to Theorem \ref{thm:primbasis2},
characterizes the basis elements $p\sub M$ of $P(H)$.
\begin{lem}\label{lem:unique}
  Suppose that $x\in P(H)$ is of the form $x=\cl M +\sum_N a\sub N\, \cl N$,
where $M$ is irreducible and each $N$ appearing in the sum is
reducible.   Then $x=p\sub M$.
\end{lem}
\begin{proof}
For any matroid $M$, we have 
$\varphi (\cl M)=\sum_{\cl N\geq\cl M}c(\cl M,\cl N)\cl N$,
where all isomorphism classes strictly greater than $\cl M$ 
appearing in the sum are reducible.  If $M$ is irreducible then
$c(\cl M, \cl M)=1$, and thus it follows that
$p\sub M= \varphi^{-1}(\cl M)$ is equal to $\cl M$ plus
a linear combination of reducible isomorphism classes of 
matroids.  Since the set \cmi\ of irreducible isomorphism classes
is linearly independent in $H$, and $\{p\sub M\co M\in\cmi\}$ is
a basis for $H$, the result follows.
\end{proof}

\begin{thm}\label{thm:primproj}
  The basis elements $p\sub M$ for $P(H)$ are given by
$$
p\sub M \= \sum_{N\in\crm}\mu\sub\ccr (M,N)\,\cl N,
$$
for all irreducible matroids $M$.
\end{thm}
\begin{proof}
The natural surjection $\pi\co C\rta H$, from the matroid-minor
coalgebra onto the matroid-minor Hopf algebra, taking a matroid $M$ to
its isomorphism class \cl M, is clearly a coalgebra map, and so
maps primitives to primitives.  Thus, if $M$ is irreducible,
it follows from
Theorem \ref{thm:primbasis2} that
$\sum_{N\in\crm}\mu\sub\ccr (M,N)\,\cl N=\pi (r\sub M)$
belongs to $P(H)$, and hence, by Lemma \ref{lem:unique}, 
$\sum_{N\in\crm}\mu\sub\ccr (M,N)\,\cl N = p\sub M$.
\end{proof}

\begin{exa}
Applying the projection $\pi$ to the primitive element $r\sub M\in P(C)$
of Example \ref{exa:lattice5}, we obtain the primitive
element $p\sub M\in P(H)$ given by Equation \ref{eq:prim5}.
\end{exa}

\begin{exa}
Applying $\pi$ to the primitive element $r\sub M$ of Example
\ref{exa:lattice7}, we obtain the primitive element
$$
\includegraphics[scale=.65,bb=94 0 450 110]
{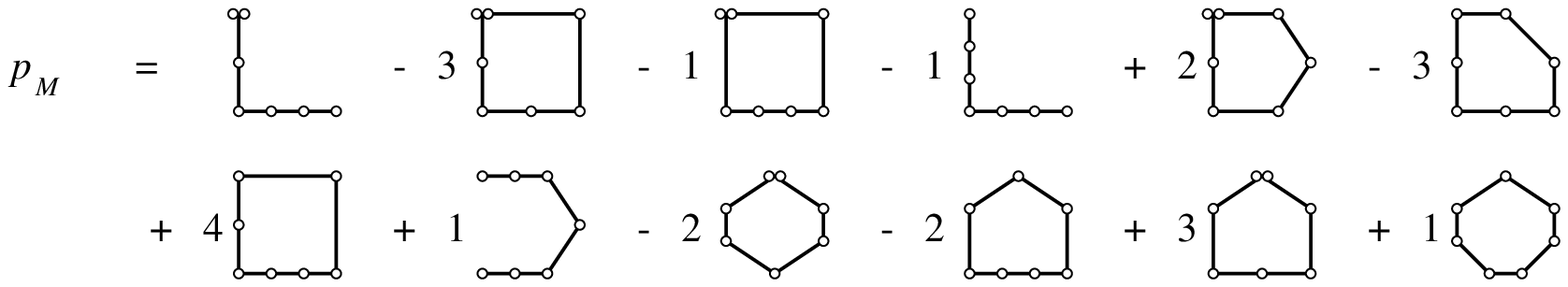}$$\\
of $H$.  This example exhibits the following curious phenomenon: the
isomorphism class of the three coatoms of the poset
${\widehat\ccr}\sub M$
does not appear in $p\sub M$, since these matroids are the only elements
of their isomorphism class belonging to
${\widehat\ccr}\sub M$
and the sum of
their M\"obius function values is zero.  We thus have a situation
in which certain matroids are needed in order to provide all the
requisite cancellations 
so that $r\sub M$ is primitive in the matroid-minor coalgebra, while
the isomorphism class of these matroids is not required in order to make
$p\sub M = \pi (r\sub M)$ primitive in the matroid-minor Hopf algebra.
\end{exa}

\section{Further work}

One line of inquiry clearly suggested by the results of the last two
sections is that of developing techniques for computing the M\"obius
functions of the posets \cwm\ and \crm.  Due to the complexity of
these posets, and the fact that so little is known about them, this is
likely to be very difficult.  A more modest, but still
worthwhile, goal would be to characterize the matroids belonging to
${\widehat\ccr}\sub M=\{N\in\crm\co\mu\sub\ccr (M,N)\neq 0\}$,
For $M$ irreducible this amounts
to identifying precisely those matroids appearing with nonzero
coefficient in the expression for the primitive element $r\sub M$ in
terms of the basis \cw\ of matroids.

A related problem, suggested by phenomenon of ``disappearing''
matroids observed in Example \ref{exa:lattice7}, is to characterize
those matroids appearing with nonzero coefficient in $r\sub M$ whose
isomorphism classes have zero coefficient in $\pi (r\sub M)$, where
$\pi\co C\rta H$ is the natural projection mapping a matroid to its
isomorphism class.
In other words, we wish to find those matroids $N$ in
${\widehat\ccr}\sub M$ such that
$\sum\mu\sub\ccr (M,N')=0$, where the sum is over all
$N'\in{\widehat\ccr}\sub M$ such that
$N'\isom N$.  Once we have a characterization of the matroids
belonging to ${\widehat\ccr}\sub M$, a solution
to this problem would provide a characterization of the isomorphism
classes that appear with nonzero coefficient in $p\sub M$.

We showed in Section \ref{sec:prim} that $r\sub M$ and $w\sub M$ are
primitive elements of $C$ whenever $M$ is irreducible (and the results
of that section imply that the converse is also true; if $M$ is
reducible, then $r\sub M$ and $w\sub M$ are not primitive).  The more
general problem of expressing the coproduct of $C$ in terms of the
bases $\{r\sub M\co M\in\cw\}$ and $\{w\sub M\co M\in\cw\}$ remains
open.  Examples we have looked at so far suggest that such expressions
will not be straightforward formulas along the lines of Equation
\ref{eq:cofreecop} but, rather, will involve some form of M\"obius
inversion and rely on the fairly subtle interplay between the free
product operation and the weak order on matroids.  We also note that
$C$ is not cofree (since it follows from Proposition \ref{pro:prodform}
that the dual algebra $C\dual$ is not free), and thus 
 the coproduct of $C$ cannot take the precise form of
Equation \ref{eq:cofreecop}.

According to Theorem \ref{thm:isom} and Equation \ref{eq:cofreecop}, 
the coproduct of $H$, expressed in terms of the basis 
$\{ p\sub M\co M\in\cm\}$ has the form
$$
\delta (p\sub M)\= \sum_{i=0}^k p\sub{M_1\frp\cdots\frp M_i}\o
p\sub{M_{i+1}\frp\cdots\frp M_k},
$$
whenever $M=M_1\frp\cdots\frp M_k$, with all $M_i$ irreducible.
Thus, in terms of this basis (and keeping in mind
unique factorization, Theorem \ref{thm:uft}), the cofreeness of $H$
becomes apparent.  We showed in Theorem \ref{thm:primproj} that the
basis element $r\sub M$ of $C$ maps to $p\sub M$ under the projection
$\pi\co C\rta H$, for irreducible $M$, thus giving, in this case, a
combinatorial interpretation for the coefficients of $p\sub M$, as
sums of M\"obius function values. The result does not hold for general
$M$, however; indeed, the coefficients of isomorphism classes
appearing in $p\sub M$ are not even necessarily integers when $M$ is
reducible (see Example \ref{exa:prim5uft}).  It would be of interest
to find a generalization of Theorem \ref{thm:primproj} that holds for
all matroids.  This would amount to determining elements of $C$, expressed
in terms of the basis $\{r\sub N\co N\in\cw\}$ (or perhaps in terms of
the basis $\{w\sub N\co N\in\cw\}$) that project to $p\sub M$, for all
matroids $M$.  A combinatorial description of such elements would
provide, in turn, a combinatorial interpretation of all coefficients
of the basis elements $p\sub M$.

Our proof of Theorem \ref{thm:primproj} is quite indirect, relying on
the uniqueness assertions of Theorem \ref{thm:primbasis2} and Lemma
\ref{lem:unique}, then using only the fact that $r\sub M$ and $p\sub
M$ are primitive, for $M$ irreducible.  A direct proof of this result
might reveal useful information about the relationship between the
posets $\crm$ and $\{\cl N\co\text{$\cl N\geq\cl M$ in \cm}\}$.

We mention finally the project of determining the precise manner in
which the matroid-minor Hopf algebra and coalgebra fit into the
framework developed by Aguiar and Sottile in \cite{agso:slr} and
\cite{agso:smr}.  A first step is to look for natural
mappings between $H$, and/or $C$, and the Hopf algebras of 
symmetric functions, permutations, and planar binary trees.

%\bibliography{discrete}

\begin{thebibliography}{10}

\bibitem{agor:hau}
Marcelo Aguiar and Rosa Orellana, \emph{The {H}opf algebra of uniform block
  permutations}, extended abstract, 2005.

\bibitem{agso:slr}
Marcelo Aguiar and Frank Sottile, \emph{Structure of the {L}oday-{R}onco {H}opf
  algebra of trees}, Journal of Algebra (2004), to appear.

\bibitem{agso:smr}
\bysame, \emph{Structure of the {M}alvenuto-{R}eutenauer {H}opf algebra of
  permutations}, Advances in Mathematics \textbf{191} (2005), no.~2, 
225--275.

\bibitem{crsc:fpm}
Henry Crapo and William Schmitt, \emph{The free product of matroids}, 
European Journal of Combinatorics \textbf{26} (2005), no.~7, 1060--1065. 
%preprint available at http://arxiv.org/abs/math.CO/0409080.

\bibitem{crsc:uft}
\bysame, \emph{A unique factorization theorem for matroids}, Journal of
  Combinatorial Theory A \textbf{112} (2005), no.~2, 222--249.
% preprint available at http://arxiv.org/abs/math.CO/0409099.

\bibitem{ge:mpp}
Ira Gessel, \emph{Multipartite {P}-partitions and inner products of skew
  {S}chur functions}, Contemporary Mathematics \textbf{34} (1984), 289--301.

\bibitem{ha:ctg}
Philip Hall, \emph{A contribution to the theory of groups of prime power
  order}, Proceedings of the London Mathematics Society \textbf{2} (1932),
  no.~36, 39--95.

\bibitem{loro:hap}
Jean-Louis Loday and Maria~O. Ronco, \emph{Hopf algebra of the planar binary
  trees}, Adv. Math. \textbf{139} (1998), 293--309.

\bibitem{loro:osc}
\bysame, \emph{On the structure of cofree {H}opf algebras}, Journal Reine u.
  Angew. Mathematik (2004), to appear.

\bibitem{ma:pcf}
Claudia Malvenuto, \emph{Produits et coproduits des fonctions
  quasisym{\'e}triques et de l'alg{\`e}bre des descents}, no.~16, Laboratoire
  de combinatoire et d'informatique math{\'e}matique (LACIM), Univ. du Quebec
  {\`a} Montr{\'e}al, 1994.

\bibitem{mare:dqf}
Claudia Malvenuto and Christophe Reutenauer, \emph{Duality between
  quasisymmetric functions and the {S}olomon descent algebra}, Journal of
  Algebra \textbf{177} (1995), no.~3, 967--982.

\bibitem{ro:fct1}
Gian-Carlo Rota, \emph{On the foundations of combinatorial theory {I}: theory
  of {M}{\"o}bius functions}, Zeitschrift f{\"u}r Wahrscheinlichkeitstheorie
  \textbf{2} (1964), 340--368.

\bibitem{sc:iha}
William Schmitt, \emph{Incidence {H}opf algebras}, Journal of Pure and Applied
  Algebra \textbf{96} (1994), 299--330.

\bibitem{we:bnm}
Dominic J.~A. Welsh, \emph{A bound for the number of matroids}, Journal of
  Combinatorial Theory \textbf{6} (1969), 313--316.

\end{thebibliography}

\providecommand{\bysame}{\leavevmode\hbox to3em{\hrulefill}\thinspace}
\providecommand{\MR}{\relax\ifhmode\unskip\space\fi MR }
% \MRhref is called by the amsart/book/proc definition of \MR.
\providecommand{\MRhref}[2]{%
  \href{http://www.ams.org/mathscinet-getitem?mr=#1}{#2}
}
\providecommand{\href}[2]{#2}

\end{document}